\documentclass[11pt]{amsart}
\usepackage{amsmath,amssymb,amsthm}
\usepackage[latin1]{inputenc}
\usepackage{version,tabularx,multicol}
\usepackage{graphicx,float}

\newcommand{\reff}[1]{(\ref{#1})}

\theoremstyle{plain}
\newtheorem{theo}{Theorem}[section]

\newtheorem{cor}[theo]{Corollary}
\newtheorem{prop}[theo]{Proposition}

\newtheorem{lem}[theo]{Lemma}

\theoremstyle{remark}
\newtheorem{rem}[theo]{Remark}

\newcommand{\cb}{{\mathcal B}}
\newcommand{\cc}{{\mathcal C}}

\newcommand{\cn}{{\mathcal N}}
\newcommand{\cm}{{\mathcal M}}
\newcommand{\cp}{{\mathcal P}}

\newcommand{\cs}{{\mathcal S}}

\newcommand{\ct}{{\mathcal T}}

\newcommand{\E}{{\mathbb E}}

\newcommand{\N}{{\mathbb N}}
\renewcommand{\P}{{\mathbb P}}

\newcommand{\R}{{\mathbb R}}

\newcommand{\rP}{{\rm P}}

\newcommand{\ind}{{\bf 1}}

\newcommand{\Supp}{{\rm Supp}\;}

\newcommand{\inv}[1]{\mathop{\frac{1}{ #1}}\nolimits}
\newcommand{\expp}[1]{\mathop {\mathrm{e}^{ #1}}}

\begin{document}

%     \cermicstitle{Fragmentation at height associated to Lévy processes}
%     \cermicsauthor[1]{Jean-François Delmas}
%     \cermicsaddress[1]{CERMICS, Ecole Nationale des
%        Ponts et Chaussées, 6-8 avenue Blaise Pascal,
%        Champs sur Marne, 77455 Marne La Vallée Cédex 2, FRANCE}
%     \makecermicstitle

\title{Fragmentation at height associated to Lévy processes}

\date{\today}

\author{Jean-François Delmas}

\address{ENPC-CERMICS, 6-8 av. Blaise Pascal,
  Champs-sur-Marne, 77455 Marne La Vallée, France.}

\email{delmas@cermics.enpc.fr}

\thanks{This research  was partially supported by 
  NSERC Discovery Grants of the Probability group at Univ. of British Columbia.}

\begin{abstract}
  We  consider the height  process of  a Lévy  process with  no negative
  jumps, and  its associated continuous tree  representation.  Using
  tools  developed  by  Duquesne and Le~Gall,  we construct a
  fragmentation process at height, which generalizes the fragmentation at
  height of stable trees given  by
  Miermont.  In this more general framework, we recover that the
  dislocation measures are the same as the dislocation measures  
  of the fragmentation at node introduced by Abraham and Delmas, up to a
  factor equal to the fragment size.  We also compute the asymptotic for
  the number of small fragments. 
\end{abstract}



\maketitle

\section{Introduction}

In  \cite{lglj:bplpep}  and \cite{lglj:bplplfss},  Le  Gall  and Le  Jan
associated to a Lévy process with  no negative jumps that does not drift
to  infinity,  $X=(X_s,  s\geq  0)$  with  Laplace  exponent  $\psi$,  a
continuous state  branching process  (CSBP) and a  Lévy continuous
random tree (CRT)  which keeps
track of the genealogy of the CSBP.  The Lévy CRT can be coded by the so
called height  process, $H=(H_s, s\geq 0)$.  Informally  $H_s$ gives the
distance (which can be understood  as the number of generations) between
the individual labeled $s$ and the root, 0, of the CRT.  We can consider
the excursion of the process $H$  above level $t>0$.  Even if $H$ is non
Markov, there  is a natural  way, we shall  recall later, to  define the
distribution $\N$ of  the excursion of $H$ above $0$  and the local time
of $H$ at level $t$ under $\N$.  Let $(\alpha_i, \beta_i)$, $i\in I$, be
the  connected component of  the open  set $\{s\in  [0,\sigma]; H_s>t\}$
($H$  is lower  semi-continuous), where  $\sigma$ is  the length  of the
excursion   of  $H$  under   $\N$.   We   denote  by   $\Lambda(t)$  the
non-increasing reordering of the sequence $(\beta_i-\alpha_i, i\in I)$,
and define the fragmentation at height: $(\Lambda(t), t\geq 0)$.

A  fragmentation process  is a  Markov  process which  describes how  an
object with  given total  mass  evolves  as it breaks  into several
fragments randomly as time passes.  Notice there may be loss of mass but
no  creation.  Those processes  have been  widely studied  in the
recent  years, see Bertoin~\cite{b:rfcp} and  references therein.  To be
more precise, the  state space of a fragmentation process  is the set of
the non-increasing sequences of masses with finite total mass
\[
\mathcal{S}^{\downarrow}=\left\{s=(s_1,s_2,\ldots); \;    s_1\ge    s_2\ge
  \cdots\ge                                         0\quad\text{and}\quad
  \Sigma(s)=\sum_{k=1}^{+\infty}s_k<+\infty\right\}.
\]
If we  denote by  $P_s$ the law  of a  $\cs^{\downarrow}$-valued process
$\Lambda       =(\Lambda       (t),t\ge       0)$      starting       at
$s=(s_1,s_2,\ldots)\in\cs^{\downarrow}$,  we say  that $\Lambda  $  is a
fragmentation  process if  it is  a Markov  process such  that $t\mapsto
\Sigma(\Lambda(t))$   is   non-increasing  and   if   it  fulfills   the
fragmentation property: the law  of $(\Lambda(t),t\ge 0)$ under $P_s$ is
the non-increasing reordering of  the fragments of independent processes
of respective  laws $P_{(s_1,0,\ldots)}$,$P_{(s_2,0,\ldots)}$, \ldots. In
other words,  each fragment  after dislocation behaves  independently of
the others,  and its evolution depends  only on its initial  mass.  As a
consequence, to describe  the law of the fragmentation  process with any
initial condition, it suffices to study the laws $P_r:=P_{(r,0,\ldots)}$
for any  $r\in (0,+\infty)$, i.e.  the law of the  fragmentation process
starting with a single mass $r$.

Theorem  \ref{theo:frag-h} states  that the  fragmentation at  height is
indeed a fragmentation process. This was already observed for the stable
case     $\psi(\lambda)=\lambda^\alpha$,    $\alpha\in     (1,2]$,    by
Bertoin~\cite{b:ssf}    ($\alpha=2$)    and   Miermont~\cite{m:sfdfstsh}
($\alpha\in (1,2)$).

A fragmentation process is said to be self-similar of index $\alpha$ if,
for  any $r>0$,  the process  $\Lambda$
under $P_r$  is distributed as the  process $(r\Lambda (r^\alpha t),t\ge
0)$  under $P_1$.  Bertoin~\cite{b:ssf}  proved that  the law  of a
self-similar   fragmentation   is  characterized   by:   the  index   of
self-similarity $\alpha$,  an erosion coefficient which  corresponds to a
 rate of  mass loss,  and a  dislocation measure  $\nu$ on
$\cs^{\downarrow}$ which describes sudden  dislocations of a fragment of
mass~1. The dislocation measure of a fragment of size $r$, $\nu_r$ is given by
$\int F(s) \nu_r(ds)= r^{\alpha} \int F(rs) \nu(ds)$. 

In the  stable cases $\psi(\lambda)=\lambda^\alpha$,  $\alpha\in (1,2]$,
the fragmentation  is self-similar with  index $-1+1/\alpha$ and  with a
zero  erosion  coefficient.   The  authors  computed in  both  cases  the
dislocation  measure and  observed it  is  the same  as the  dislocation
measure   associated  to   the   fragmentation  at   ``nodes''  of   the
corresponding CRT  see \cite{ap:sac} and  \cite{b:fpcbm} for $\alpha=2$,
and \cite{m:sfdfstsn} for $\alpha \in (1,2)$.

For a  general sub-critical  or critical CRT,  there is no  more scaling
property, and the dislocation measure, which describes how a fragment of
size  $r>0$ is  cut in  smaller  pieces, can't  be expressed  as a  nice
function  of  the  dislocation measure  of  a  fragment  of size  1.  In
\cite{ad:falp},  the  authors give  the  family  of dislocation  measures
$(\nu_r, r>0)$ for the fragmentation at node of a general sub-critical or
critical    CRT.     We     set    $\nu^*_r=r^{-1}    \nu_r$.    Theorem
\ref{theo:intensite} state that $(\nu_r^*, r>0)$  is the family  of
dislocation  measures for the fragmentation at height, see also Remark
\ref{rem:frag-int}. Intuitively $\nu_r^*$ describes the way a mass $r$
breaks in smaller pieces. 

We also compute the asymptotic of  the number of small fragments at time
$t$  (see   \cite{b:smssf}  and   \cite{h:rfdsf}  for  results   in  the
self-similar case). With a suitable renormalization, it converges to the
local  time  of  the  height  process  at  level  $t$,  see  Proposition
\ref{prop:smallN}  and Corollary  \ref{cor:smallN} for  the  stable case
$\alpha\in (1,2)$.  We  also characterize the law of  this local time at
level $t$ under $\rP_r$ in Lemma \ref{lem:L_s=1}.

The paper is organized as follows. In Section \ref{sec:not}, we recall
the definition and properties of the height and exploration
processes. In the very short Section \ref{sec:fh} we state and prove the
fragmentation at height is indeed a fragmentation process. The number of
small fragments is studied in Section \ref{sec:smallN}. And the
dislocation measures are computed in Section \ref{sec:dm}.

\section{Notations}
 \label{sec:not}
 We denote by $\cb_+(\R_+)$ the set of measurable non-negative functions
 defined  on $\R_+$.   Let  $\cm_f$ be  the  set of  finite measures  on
 $\R_+$,  endowed with  the topology  of weak  convergence.  For $\mu\in
 \cm_f$, we set $H^\mu=\sup \Supp (\mu)\in [0, \infty ]$ the supremum of
 its closed  support.  For $f\in \cb_+(\R_+)$,  we write
 $\langle \mu, f \rangle$ for $\int_{\R_+} f(x)\; \mu(dx)$.

Let $\psi$ denote the Laplace
exponent            of            $X$:           $\E\left[\expp{-\lambda
    X_t}\right]=\expp{t\psi(\lambda)}$,  $\lambda>0$.   We shall  assume
there is no Brownian part, so that
\[
\psi(\lambda)=\alpha_0\lambda+\int_{(0,+\infty)}\pi(d\ell)
\left[\expp{-\lambda\ell}-1+\lambda\ell\right],  
\]
with  $\alpha_0\ge  0$  and  the   Lévy  measure  $\pi$  is  a  positive
$\sigma$-finite measure  on $(0,+\infty)$ such  that $\int_{(0,+\infty)}
(\ell\wedge \ell^2)\pi(d\ell)<\infty$.  Following \cite{dlg:rtlpsbp}, we
shall also assume that $X$  is of infinite variation a.s.  which implies
that  $\int_{(0,1)}\ell\pi(d\ell)=\infty$.  Notice those  hypothesis are
fulfilled in the stable case: $\psi(\lambda)=\lambda^\alpha$, $\alpha\in
(1,2)$.      For      $\lambda\geq     1/\varepsilon>0$,     we     have
$\expp{-\lambda\ell}-1+ \lambda\ell\geq \frac{1}{2}\lambda \ell \ind_{\{
  \ell\geq   2\varepsilon\}}$,    which   implies   that   $\lambda^{-1}
\psi(\lambda)  \geq  \alpha_0+  \int_{(2\varepsilon,\infty  )}  \ell  \;
\pi(d\ell)$. We deduce that
\begin{equation}
   \label{eq:psi/l}
\lim_{\lambda
\rightarrow\infty } \frac{\lambda}{\psi(\lambda)} =0.
\end{equation}

The   so-called  exploration  process   $\rho=(\rho_t,t\ge  0)$   is  an
$\cm_f$-valued càd-làg  Markov process.  The height process  at time $t$
is defined as  the supremum of the closed support  of $\rho_t$ (with the
convention  that  $H_t=0$ if  $\rho_t=0$).   Informally, $H_t$  gives the  distance (which  can be
understood as the number  of generations) between the individual labeled
$t$ and the root, 0, of the CRT.  In some sense $\rho_t(dv)$ records the
``number'' of brothers, with labels  larger than $t$, of the ancestor of
$t$ at generation $v$.

We recall the definition and properties of the exploration process which
are    given    in    \cite{lglj:bplpep},    \cite{lglj:bplplfss}    and
\cite{dlg:rtlpsbp}. The results of  this section are mainly extracted from
\cite{dlg:rtlpsbp}.
 
Let $I=(I_t,t\ge 0)$ be the infimum process of $X$, $I_t=\inf_{0\le s\le
  t}X_s$.  We will also consider for  every $0\le s\le t$ the infimum of
$X$ over $[s,t]$:
\[
I_t^s=\inf_{s\le r\le t}X_r.
\]
There exists a sequence $(\varepsilon_n,n\in \N^*)$ of positive real
numbers decreasing to 0 s.t. 
\[
\tilde H_t= \lim_{k\rightarrow\infty } \inv{\varepsilon_k} \int_0^t
\ind_{\{X_s<I^s_t+\varepsilon_k\}}\; ds
\]
exists and is finite a.s. for all $t\geq 0$. 

The point 0  is regular for the Markov process $X-I$,  $-I$ is the local
time  of $X-I$  at 0  and the  right continuous  inverse of  $-I$  is a
subordinator  with   Laplace exponent $\phi$,  the   inverse  of  $\psi$:
$\psi(\phi(x))=\phi(\psi(x))=x$  (see  \cite{b:pl},  chap. VII).  Notice
this subordinator  has no drift  thanks to \reff{eq:psi/l}.  Let $\pi_*$
denote the corresponding Lévy measure.

Let $\N$ be the associated excursion measure of the process $X-I$ out of
0, and $\sigma=\inf\{t>0; X_t-I_t=0\}$ be the length of the excursion of
$X-I$ under  $\N$. Under  $\N$, $X_0=I_0=0$. We  shall use  (see Section
3.2.2. in \cite{dlg:rtlpsbp}) that
\begin{equation}
   \label{eq:exp-sigma}
\N[1-\expp{-\lambda \sigma}]=\phi(\lambda).
\end{equation}
In particular $\sigma$  is distributed under $\N$ according  to the Lévy
measure $\pi_*$.

From section  1.2 in \cite{dlg:rtlpsbp},  there exists a  $\cm_f$ valued
process, $\rho^0=(\rho^0_t,  t\geq 0)$, called  the exploration process,
such that :
\begin{itemize}
\item A.s., for every $t\geq 0$, we have 
$\langle \rho_t^0,1\rangle =X_t-I_t$, and the process $\rho^0$ is càd-làg. 
\item The process $(H_s^0=H^{\rho^0_s}, s\geq 0)$ taking values in
  $[0,\infty ]$ is lower semi-continuous. 
\item   For  each   $t\geq   0$,  a.s.   $H^0_t=   \tilde  H_t$.
\item For every $f\in \cb_+(\R_+) $,
\[
\langle \rho^0_t,f\rangle =\int_{[0,t]} f(H^0_s)\; d_sI_t^s,
\]
or equivalently, with $\delta_x$ being the Dirac mass at $x$, 
\[
\rho^0_t(dr)=\sum_{\stackrel{0<s\le t}
  {X_{s-}<I_t^s}}(I_t^s-X_{s-})\delta_{H_s^0}(dr). 
\]
\end{itemize}

In the definition  of the exploration process, as $X$  starts from 0, we
have  $\rho_0=0$ a.s.  To  get the  Markov property  of $\rho$,  we must
define  the  process  $\rho$  started  at any  initial  measure  $\mu\in
\cm_f$.  For  $a\in [0,  \langle \mu,1\rangle ]  $, we  define the
erased measure $k_a\mu$ by
\[
k_a\mu([0,r])=\mu([0,r])\wedge (\langle \mu,1\rangle -a), \quad
\text{for $r\geq 0$}. 
\]

If $a> \langle  \mu,1\rangle $, we set $k_a\mu=0$.   In other words, the
measure $k_a\mu$ is the measure $\mu$ erased by a mass $a$ backward from
$H^\mu$.

For $\nu,\mu \in \cm_f $, and $\mu$ with compact support, we define
the concatenation $[\mu,\nu]\in \cm_f  $ of the two measures by:
\[
\bigl\langle [\mu,\nu],f\bigr\rangle =\bigl\langle \mu,f\bigr\rangle
+\bigl\langle \nu,f(H^\mu+\cdot)\bigr\rangle , 
\quad f\in \cb_+(\R_+).
\]
Eventually, we  set for every $\mu\in \cm_f $ 
and every $t>0$, 
\[
\rho_t=\bigl[k_{-I_t}\mu,\rho_t^0].
\]
We say  that $\rho=(\rho_t, t\geq 0)$  is the process  $\rho$ started at
$\rho_0=\mu$, and  write $\P_\mu$ for its law.  We set $H_t=H^{\rho_t}$.
The process $\rho$ is càd-làg and strong Markov.

\subsection{Poisson decomposition}
\label{sec:pois_d}
Let $\P^*_\mu$ denote the law of $\rho$ started at $\mu$ and killed when
it  reaches $0$.   We decompose  the path  of $\rho$  under $\P^*_{\mu}$
according to excursions  of the total mass of  $\rho$ above its minimum,
see   Section   4.2.3  in   \cite{dlg:rtlpsbp}.    More  precisely   let
$(\alpha_i,\beta_i), i\in  I$ be the excursion intervals  of the process
$X-I$ away from  $0$ under $\P^*_{\mu}$.  For every  $i\in I$, we define
$h_i=H_{\alpha_i}$     and     $\rho^i$     by     the     formula     $
\rho_t^i=\rho_{(\alpha_i+t)\wedge    \beta_i}   ^0$    or   equivalently
$[k_{-I_{\alpha_i}} \mu, \rho^i_t]=\rho_{(\alpha_i+t)\wedge \beta_i}$.
We recall  Lemma  4.2.4. of \cite{dlg:rtlpsbp}.
\begin{lem}
\label{lem:dlg-decomp}
   Let $\mu\in \cm_f $. The point
   measure $\displaystyle \sum_{i\in I} \delta_{(h_i,\rho^i)}$ is under
    $\P^*_{\mu}$ a Poisson point measure with intensity $\mu(dr) \N[d\rho]$.
\end{lem}

Let $(L_s^t,s\geq  0, t>  0)$ be  the local time  of the  height process
under     $\N$    at     level    $t>0$     at    time     $s\geq    0$:
$L^t_s=\lim_{\varepsilon\rightarrow    0}   \inv{\varepsilon}   \int_0^s
\ind_{\{t<H_r\leq t+\varepsilon\}} dr$ in $L^1$-norm.  The local time is
jointly measurable in $(s,t)$, non-decreasing and continuous in $s$ (see
Proposition 1.3.3 in \cite{dlg:rtlpsbp}).

Consider the  right-continuous inverse of  the time spent by  the height
process  under  level  $t$:  $\tilde  \tau_s^t=\inf\{r\geq  0;  \int_0^r
\ind_{\{ H_v\leq  t\}} \;dv >s\}$, and set  $\tilde \rho=(\tilde \rho_s,
s\geq   0)$  where   $\tilde  \rho_s=\rho_{\tilde   \tau_s^t}$.   Recall
$L^t_\sigma$  is measurable  with  respect to  $\tilde  \rho$ thanks  to
Proposition 1.3.3 in \cite{dlg:rtlpsbp}.

Notice the set  $\{s\in [0,\sigma]; H_s>t\}$ is open  since $H$ is lower
semi-continuous.  Let $(\alpha_i, \beta_i)$, $i\in I_t$ be the excursion
of  $H$ (or  $\rho$) above  level $t$.   Notice that  $I_t=\emptyset$ if
$\sup \{  H_s, s\in  [0,\sigma] \}<t$. If  $I_t\neq \emptyset$,  we define
$\rho^i$  such   that  $[\rho_{\alpha_i},  \rho^i_t]=\rho_{(\alpha_i  +t
  )\wedge \beta_i}$, and $\sigma^i=\beta_i-\alpha_i$ the duration of the
excursion  $\rho^i$.   By Proposition  1.3.1  in \cite{dlg:rtlpsbp}  and
standard excursion theory, we have the next Lemma.

\begin{lem}
   \label{lem:pm>t}
   Under $\N$,  the point measure $\sum_{i\in  I_t} \delta_{\rho^i}$ is,
   conditionally on $L_\sigma^t $ (or  on $\tilde \rho$) a Poisson point
   measure with intensity $L^t_\sigma \N[d\rho]$.
\end{lem}

\subsection{The dual process and representation formula}
\label{sec:dual}

We  shall need the  $\cm_f $-valued process  $\eta=(\eta_t,t\ge 0)$
 defined by
\[
\eta_t(dr)=\sum_{\stackrel{0<s\le t}{X_{s-}<I_t^s}}(X_s-I_t^s)\delta
_{H_s}(dr).
\]
The process $\eta$ is the dual process of $\rho$ under $\N$ (see
Corollary 3.1.6 in \cite{dlg:rtlpsbp}). Let $\Delta_s=X_s -X_{s-}$,
$s>0$, be the jumps of $X$. We write  (recall
$\Delta_s= X_s-X_{s-}$)
\[
\kappa_t(dr)=\rho_t(dr)+\eta_t(dr)=\sum_{\stackrel{0<s\le
    t}{X_{s-}<I_s^t}}\Delta_s \delta _{H_s}(dr).
\]

We recall the Poisson representation of $(\rho,\eta)$ under $\N$. Let
$\mathcal{N}(dx\,   d\ell\,  du)$   be  a   Poisson  point   measure  on
$[0,+\infty)^3$ with intensity
$$dx\,\ell\pi(d\ell)\ind_{[0,1]}(u)du.$$
For every $a>0$, let us denote by $\mathbb{M}_a$ the law of the pair
$(\mu_a,\nu_a)$ of finite measures on $\R_+$ defined by:  for $f\in \cb_+(\R_+)$ 
\begin{align*}
\langle \mu_a,f\rangle  & =\int\mathcal{N}(dx\, d\ell\, du)\ind_{[0,a]}(x)u\ell f(x),\\
\langle \nu_a,f\rangle  & =\int\mathcal{N}(dx\, d\ell\, du)\ind_{[0,a]}(x)\ell(1-u)f(x).
\end{align*}
We eventually set $\mathbb{M}=\int_0^{+\infty}da\, \expp{-\alpha_0
  a}\mathbb{M}_a$. 
For every non-negative measurable function $F$ on $\cm_f ^2$, we have 
\begin{equation}
   \label{eq:poisson-rh}
\N\left[\int_0^\sigma F(\rho_t, \eta_t) \; dt \right]=\int\mathbb{M}(d\mu\,
    d\nu)F (\mu, \nu) ,
\end{equation}
where $\sigma=\inf\{s>0; \rho_s=0\}$ denotes the length of the
    excursion.  
Let $t>0$.  For every non-negative measurable function
$F$ on $\cm_f $, we have
\begin{equation}
   \label{eq:poisson-rl}
\N\left[\int_0^\sigma F(\rho_s) \; dL^t_s \right]=\expp{-\alpha_0 t} 
\int \mathbb{M}_t(d\mu\,
    d\nu)F (\mu).
\end{equation}

\section{The fragmentation at height}
\label{sec:fh}

We   keep   notations   introduced   for   Lemma   \ref{lem:pm>t} in
section \ref{sec:pois_d}.   We define at  time $t$ 
the     fragmentation      process     at     height,     $\Lambda(t)\in
\mathcal{S}^{\downarrow}$, as the sequence $(\sigma^i, i\in I_t)$ ranked
in non-increasing order  (if $I_t$ is empty or  finite, this sequence is
completed  by  zeroes).   If  needed,  we  write  $\Lambda^\rho(t)$  for
$\Lambda(t)$  to stress  that the  fragmentation process  is  built from
$\rho$.

Let  $\pi_*$ be the distribution of $\sigma$ under $\N$. 
By decomposing the measure $\N$ w.r.t.  the distribution of $\sigma$, we
get that  $\N[d\rho]=\int _{(0,\infty )} \pi_*  (dr) \N_r[d\rho]$, where
$(\N_r,  r \in  (0,\infty  ))$  is a  measurable  family of  probability
measures on the  set of excursions of the  exploration process such that
$\N_r[\sigma=r]=1$ for $\pi_*(dr)$-a.e. $r>0$. From standard excursion
theory, we have the next Lemma. 

\begin{lem}
   \label{lem:rho_i>t}
Let $t>0$. The random variables  $(\rho_i, i\in I_t)$ are, conditionally on the
$\sigma$-field generated by 
$(H_s \wedge 
t, s\in [0,\sigma])$ and $((\alpha_i, \beta_i), i \in I_t)$, independent and
$\rho^i$ is distributed according to $\N_{\sigma^i}[d\rho]$. 
\end{lem}

For  $\pi_*(dr)$-a.e., let  $\rP_r$ denote  the law  of  $\Lambda$ under
$\N_r$, and let $\rP_0$ be the law of the constant process equal to $(0,
\ldots)\in \mathcal{S}^{\downarrow}$.

\begin{theo}
\label{theo:frag-h}
  For $\pi_*(dr)$-a.e.,  under $\rP_r$, the  process $\Lambda
$    is    a    $\mathcal{S}^\downarrow$-valued
  fragmentation process.   More precisely,  $\Lambda$ is Markov  and
  satisfy the fragmentation property: the
  law  under   $\rP_r$  of   the  process  $(\Lambda   (t+t'),t'\ge  0)$
  conditionally on $\Lambda (t)=(\Lambda _i, i\in \N^*)$ is given
  by the  decreasing reordering of independent  processes of respective
  law $\rP_{\Lambda _i}$, $i\in  \N^*$. 
\end{theo}
\begin{proof}
Notice that $(\Lambda(s), s\in [0,t])$  is measurable w.r.t. the
$\sigma$-field generated by $(H_u\wedge t, u\in [0,\sigma])$ and
$((\alpha_i, \beta_i), i\in I_t)$. Using notations of Lemma
\ref{lem:rho_i>t},  for $t'>0$, $\Lambda^\rho(t+t')$ is the
non-increasing reordering
of $(\Lambda^{\rho^i}(t'), i\in I_t)$. The Markov property and the
fragmentation property are consequences of 
Lemma \ref{lem:rho_i>t}.
\end{proof}

\section{Number of small fragments}
\label{sec:smallN}

For the  fragmentation at height, it  is easy to give  the asymptotic of
the number of  small fragments.  We keep notations  introduced for Lemma
\ref{lem:pm>t}.   Let $N_\varepsilon(t)$  be the  number of  fragments at
time $t$ of size bigger or equal to $\varepsilon>0$ and
$M_\varepsilon(t)$ the total mass of the fragments less or equal than
$\varepsilon$:
\[
N_\varepsilon(t)=\sum_{i\in I_t} \ind_{\{ \sigma^i\geq
  \varepsilon\}}\quad\text{and}\quad  M_\varepsilon(t)=\sum_{i\in I_t}
\sigma^i \ind_{\{ \sigma^i\leq
  \varepsilon\}}.
\]
For $t>0$, we write $\bar \pi_*(t)=\pi_*((t,\infty ))$ and
$\varphi(t)=\int_{(0,t)} r\pi_*(dr)$. 

\begin{lem}
\label{lem:cv-f/e}
We have $\lim_{\varepsilon \rightarrow 0} \bar \pi_*(\varepsilon)=\infty
$ and $\lim_{\varepsilon \rightarrow 0} \inv{\varepsilon} 
\varphi(\varepsilon)=\infty $. 
\end{lem}

\begin{proof}
  We deduce from \reff{eq:psi/l} that 
\begin{equation}
   \label{eq:phi}
\lim_{\lambda\rightarrow \infty }
  \phi(\lambda)=\infty. 
\end{equation} 
Notice   $\phi(\lambda)=   \int_{(0,\infty   )}  (1-\expp{-\lambda   r})
\pi_*(dr)$ to obtain the first part of the Lemma.

The     second     limit      is     more     involved.      We     have
$\lambda\psi'(\lambda)=\alpha_0 \lambda  + \int_{(0,\infty )}  \lambda r
(1-\expp{-\lambda r}) \pi(dr)$. Since for all $x\geq 0$,
\[
\expp{-x} -1 +x \leq  x(1-\expp{-x}) \leq  2\left(\expp{-x} -1 +x
\right),
\]
we deduce that $ \psi(\lambda)\leq \lambda \psi'(\lambda) \leq  2
\psi(\lambda)$. And we get that for $x>0$, 
\begin{equation}
   \label{eq:majo-phi}
\inv{2} \phi(x) \leq  x\phi'(x) \leq \phi(x).
\end{equation}
The  next  part  of the  proof is
  inspired  by a  Theorem  of Haan  and  Stadtmüller (see  \cite{bgt:rv}
  p.118). We have 
\[
\phi'(\lambda)=\int_{(0,\infty )} \expp{-\lambda r } r \pi_*(dr)=
\int_{(0,\infty )} \expp{-u} \varphi(u/\lambda)\; du.
\]
Since the function $\varphi$ is non-decreasing, we have for $z>0$
\begin{equation}
   \label{eq:majo_f}
\phi'(\lambda)\geq \int_{(z,\infty )} \expp{-u} \varphi(u/\lambda)\;
du\geq  \expp{-z} \varphi(z/\lambda). 
\end{equation}
We have for any $x>0$, 
\begin{align*}
   \inv{2\lambda} \phi(\lambda)
&\leq  \phi'(\lambda)\\
&\leq \int_{(0,x )} \expp{-u} \varphi(u/\lambda)\; du+ \int_{(x,\infty
  )} \expp{-u} \varphi(u/\lambda)\; du \\
&\leq \varphi(x/\lambda) ( 1- \expp{-x}) +  \int_{(x,\infty
  )} \expp{-u} \expp{u/2} \phi'(\lambda/2)\; du \\
&\leq \varphi(x/\lambda) +  \frac{4}{\lambda} \phi(\lambda/2)\expp{-x/2},
\end{align*}
where we used \reff{eq:majo-phi} for the first inequality, the fact
that $\varphi$ is non-decreasing for the first term and \reff{eq:majo_f}
for the second term of the third inequality, and \reff{eq:majo-phi}
again for the last. With $x=6\ln(2)$, we get
\[
\varphi(x/\lambda) \geq \inv{2\lambda} (\phi(\lambda)-\phi(\lambda/2))
=\inv{2\lambda}\int_{[\lambda/2,\lambda]} \phi'(u)\; du
\geq  \frac{\phi'(\lambda)}{4},
\]
where we used that $\phi'$ is non-increasing for the last
inequality. Using \reff{eq:majo-phi}, we deduce that 
\[
\lambda\varphi(1/\lambda)\geq \frac{\lambda}{4} \phi'(\lambda/x) \geq
\frac{x}{8}\phi(\lambda/x).
\]
The last part of the Lemma is then a consequence of \reff{eq:phi}. 
\end{proof}

\begin{prop}
\label{prop:smallN}
Let   $t>0$.     We   have   that,    conditionally   on   $L_\sigma^t$,
$N_\varepsilon(t)$   is   a    Poisson   random   variable   with   mean
$\bar\pi_*(\varepsilon)L_\sigma^t$.  Furthermore, there  exists a
sequence  of positive  numbers, $(\varepsilon_n,  n\geq  1)$, decreasing
towards 0, such that, for all $t>0$, we have $\N$-a.e.
\[
\lim_{n\rightarrow \infty }
\frac{N_{\varepsilon_n}(t)}{\bar \pi_*(\varepsilon_n)} 
=\lim_{n\rightarrow \infty }
\frac{M_{\varepsilon_n}(t)}{\varphi(\varepsilon_n)} 
=L_\sigma^t.
\]
\end{prop}
We can replace  the sequence $(\varepsilon_n, n\geq 1)$  by any sequence
in some  particular case  (see for example  Corollary \ref{cor:smallN}).
Those results extend \cite{h:rfdsf} for the fragmentation at height.

\begin{proof}
  Recall $\N[\sigma \geq  \varepsilon]= \bar \pi_*(\varepsilon)$.  Since
  conditionally  on $L^t_\sigma$,  $\sum_{i\in I_t}  \delta_{\rho^i}$ is
  under   $\N$   a   Poisson   point  measure   with   intensity   $\bar
  \pi_*(\varepsilon)L_\sigma^t$  (see Lemma  \ref{lem:pm>t}),  we deduce
  that  conditionally on $L_\sigma^t$,  $N_\varepsilon(t)$ is  a Poisson
  random variable with mean $\bar \pi_*(\varepsilon)L_\sigma^t$.

Since $\lim_{\varepsilon\downarrow 0} \bar\pi_*(\varepsilon)=\infty
$, there exists a sequence $(\varepsilon_n, n\geq 1)$ decreasing
towards 0 such that $\sum_{n\geq 1} 1/ \bar \pi_*(\varepsilon_n)<\infty
$. Since $\displaystyle
\N\left[\left(\frac{N_\varepsilon(t)}{\bar \pi_*(\varepsilon)} -
    L_\sigma^t  \right)^2 \Big| L_\sigma^t \right]=
\frac{L_\sigma^t}{ 
 \bar \pi_*(\varepsilon)}$, the series $\displaystyle \sum_{n\geq 1}
\left(\frac{N_{ \varepsilon_n}(t)}{\bar \pi_*(\varepsilon_n)} -
    L_\sigma^t  \right)^2$ 
is $\N$-a.e. finite. This implies that $\N$-a.e. $\displaystyle  \lim_{n\rightarrow \infty }
\frac{N_{\varepsilon_n}(t)}{\bar \pi_*(\varepsilon_n)} 
=L_\sigma^t$.

Since conditionally  on $L^t_\sigma$,  $\sum_{i\in I_t}  \delta_{\rho^i}$ is
  under   $\N$   a   Poisson   point  measure   with   intensity   $\bar
  \pi_*(\varepsilon)L_\sigma^t$  (see Lemma  \ref{lem:pm>t}),  we deduce
  that for $\lambda>0$, $\varepsilon>0$,
\[
\N\left[\expp{- \lambda M_\varepsilon(t)}
    \Big| L^t_\sigma\right]=\expp{- L_\sigma^t \int_{(0,\varepsilon)}
    \pi_*(dr)\; (1-\expp{-\lambda r})}.
\]
Poisson point measure properties yield 
\begin{align*}
\N\left[M_\varepsilon(t)
    | L^t_\sigma\right]&= L_\sigma^t \int_{(0,\varepsilon )}
   r \pi_*(dr)=\varphi(\varepsilon) L^t_\sigma\\
\N\left[M_\varepsilon(t)^2
    | L^t_\sigma\right]&= \N\left[M_\varepsilon(t)
    | L^t_\sigma\right]^2+ L_\sigma^t \int_{(0,\varepsilon )}
    r^2 \pi_*(dr).
\end{align*}
We deduce that 
\[
 \N\left[\left(\frac{M_\varepsilon(t)}{\varphi(\varepsilon)}
     -L^t_\sigma \right)^2\Big | L^t_\sigma \right]
=  L_\sigma^t \frac{\int_{(0,\varepsilon )}
    r^2 \pi_*(dr)}{\varphi(\varepsilon)^2}
\leq  L_\sigma^t \frac{\varepsilon \int_{(0,\varepsilon )}
    r \pi_*(dr)}{\varphi(\varepsilon)^2} =
  \frac{\varepsilon}{\varphi(\varepsilon)} L^t_\sigma.
\]
As $\lim_{\varepsilon\rightarrow 0} \varepsilon/\varphi(\varepsilon)=0$
(see Lemma \ref{lem:cv-f/e}), we can use similar arguments as those used
for $N_\varepsilon(t)$ to end the proof. Eventually, notice one can
choose the sequence $(\varepsilon_n, n\geq 1)$ such that the two limits
in the Proposition hold simultaneously. 
\end{proof}

\begin{rem}
As $\lim_{\varepsilon\rightarrow 0} \varepsilon/\varphi(\varepsilon)=0$
(see Lemma \ref{lem:cv-f/e}), we have 
\[
\lim_{\varepsilon\rightarrow 0} \int_{(0,\infty )}
    \pi_*(dr)\; (1-\expp{-\lambda r/\varphi(\varepsilon)}) = \lambda.
\]
This implies that from the above proof
\[
\lim_{\varepsilon\rightarrow 0} \N\left[\expp{- \lambda
    M_\varepsilon(t)/\varphi(\varepsilon) }
    \Big| L^t_\sigma\right]=\expp{- L_\sigma^t
\displaystyle     \lim_{\varepsilon\rightarrow 0} \int_{(0,\varepsilon)}
    \pi_*(dr)\; (1-\expp{-\lambda r/\varphi(\varepsilon)})} =\expp{-
    \lambda L_\sigma^t}.  
\]
This   implies      that      conditionally      on      $L^t_\sigma$,
$M_\varepsilon(t)/\varphi(\varepsilon)$  converges   in  probability  to
$L^t_\sigma$  as  $\varepsilon$  goes   down  to  0.  Notice  also  that
conditionally        on       $L^t_\sigma$,       $N_\varepsilon(t)/\bar
\pi_*(\varepsilon)$   converges  in   probability  to   $L^t_\sigma$  as
$\varepsilon$ goes down to 0.
\end{rem}

We  consider the  stable case  $\psi(\lambda)=\lambda^\alpha$.   We have
$\pi_*(dr)= (\alpha  \Gamma(1-\alpha^{-1}))^{-1} r^{-1-1/\alpha} \; dr$,
$\bar           \pi_*(\varepsilon)=           \Gamma(1-\alpha^{-1})^{-1}
\varepsilon^{-1/\alpha}$,   and  $\varphi(\varepsilon)=  \Big((\alpha-1)
\Gamma(1-\alpha^{-1})\Big)^{-1}  \varepsilon^{1-\alpha^{-1}}$.  There is
a version of $(\N_r, r>0)$ such that for all $r>0$ we have $\N_r[F((X_s,
s\in  [0,r]))]=\N_1[F((r^{1/\alpha}  X_{s/r} ,  s\in  [0,r]))]$ for  any
non-negative  measurable function  $F$  defined on  the  set of  càd-làg
paths. For the next Corollary, see also general results from  \cite{h:rfdsf}.

\begin{cor}
\label{cor:smallN}
  Let $\psi(\lambda)=\lambda^\alpha$, for  $\alpha\in (1,2)$. Let $t>0$.
  We  have, under $\N$  or $\N_1$,  that conditionally  on $L_\sigma^t$,
  $N_\varepsilon(t)$ is a Poisson  random variable with mean
  $\Gamma(1-\alpha^{-1} )^{-1} \varepsilon^{-1/\alpha} L_\sigma^t$.  Furthermore,
 for all $t>0$, we have $\N$-a.e. or $\N_1$-a.s.
\[
\lim_{\varepsilon \rightarrow 0} \varepsilon^{1/\alpha} 
N_{\varepsilon}(t)
=\lim_{\varepsilon \rightarrow 0} (\alpha-1) \varepsilon^{\alpha^{-1} -1} 
M_{\varepsilon}(t)
= \frac{L_\sigma^t} {\Gamma(1-\alpha^{-1})}.
\]
\end{cor}

\begin{proof}
   It is enough to check that we can replace the sequence
   $(\varepsilon_n, n\geq 1)$ in Proposition \ref{prop:smallN} by any
   sequences. Notice that $\varepsilon_n=n^{-2\alpha}$, $n\geq 1$,
   satisfies $\sum_{n\geq 1} 1/\bar\pi_*(\varepsilon_n)<\infty
   $. From the proof of Proposition \ref{prop:smallN}, we get that
   $\N$-a.e. or $\N_1$-a.s., $\lim_{n \rightarrow\infty } n^{-2}
   N_{n^{-2\alpha}}(t) =L^t_\sigma/ \Gamma(1-\alpha^{-1})$. Since
   $N_\varepsilon(t)$ is a non-increasing function of $\varepsilon$, we
   get that for any $\varepsilon\in [(n+1)^{-2\alpha}, n^{-2\alpha}]$,
   we have 
\[
\frac{n^2}{(n+1)^2} n^{-2}
   N_{n^{-2\alpha}}(t) \leq  \varepsilon^{1/\alpha} N_\varepsilon(t)
\leq  \frac{(n+1)^2}{n^2} (n+1)^{-2}
   N_{(n+1)^{-2\alpha}}(t) .
\]
   Hence we deduce that  $\N$-a.e. or $\N_1$-a.s., $\lim_{\varepsilon
     \rightarrow\ 0 } \varepsilon^{1/\alpha} 
   N_{\varepsilon}(t) =L^t_\sigma/ \Gamma(1-\alpha^{-1})$. The proof for
   $M_\varepsilon(t)$ is similar, as $M_\varepsilon(t)$ is a
   non-decreasing function of $\varepsilon$. 
\end{proof}

The next Lemma characterizes the law of $L^t_\sigma$ under
$\N_r$. (Recall $\int F(r, \rho)  \;  \pi_*(dr)\N_r[d\rho]=\N[F(s,\rho)]$.) 

\begin{lem}
   \label{lem:L_s=1}
Let $\lambda\geq 0$ and $\gamma\geq 0$. 
   Let   $w(t)=\N\left[\expp{-\lambda   \sigma}  \left(1   -\expp{-\gamma
         L_\sigma^t}\right) \right]$,  for $t>0$ and
   $w(0)=\gamma$.  Then $w$ belongs to  
     $\cc^1(\R_+)$, is non-increasing, such that $\lim_{t\rightarrow\infty }
     w(t)=0$ and  solves 
 \[
w'(t)=\lambda  -  \psi(\phi(\lambda)  +
     w(t)), \quad t>0.
\]
\end{lem}

\begin{rem}
  For the  $\alpha$-stable case,  we can characterize  the  law of
  $L^t_\sigma$  under $\N_1$.  Scaling properties  yield  the processes
  $r^{1/\alpha}  X_{s/r}$,  $r^{1-\alpha^{-1}}  H_{s/r}$,  $r^{1/\alpha}
  \rho_{s/r}   (d   (u   \;r^{1-  \alpha^{-1}}))$ and   $(r^{1/\alpha}
  L_{s/r}^{t\;r^{-1+\alpha^{-1}} }\!\!  , t\geq  0)$ are distributed  as the
  processes $X_s$, $H_s$, $\rho_s(du)$ and $(L_s^t, t\geq 0)$. We deduce
  from the definition of $w$ in Lemma \ref{lem:L_s=1} that
\[
\inv{\alpha \Gamma (1- \alpha^{-1})} \int_{(0,\infty )} dr\;
r^{-1-1/\alpha} \expp{-\lambda r} \left[1 - \N_1\left[\expp{- \gamma
      r^{1/\alpha} L^{tr^{(-1+\alpha^{-1})}}_1}\right]\right]=w(t).
\]
And this is enough to characterize the law of $L^t_1$ under $\N_1$, for
all $t>0$. 
\end{rem}

\begin{proof}[Proof of Lemma \ref{lem:L_s=1}.]
We first consider for $\lambda>0$,  $\gamma>0$ and $t>0$ 
\[
v(t)=\N\left[1 -\expp{-\lambda \sigma -\gamma L^t_\sigma}\right]   .
\]

Notice that $\sigma=\sum_{i\in I_t} \sigma^i + \tilde \sigma$, where
$\tilde \sigma=\int_0^\sigma \ind_{\{H_s\leq t\}} \; ds $ is the duration of the excursion of $\tilde \rho$
(defined in Section \ref{sec:pois_d}) under $\N$. 
From Lemma \ref{lem:pm>t}, we get that by conditioning with respect to
$\tilde \rho$,
\[
v(t)=\N\left[1 -\expp{-\lambda \int_0^\sigma \ind_{\{H_s\leq t\}} \; ds
    - L^t_\sigma(\gamma+ \N[1-\expp{-\lambda \sigma}] ) } \right]   .
\]
We consider the additive functional of $\rho$ given by
\[
dA_s=\lambda \ind_{\{H_s\leq t\}} \; ds + (\gamma+
\phi(\lambda))dL_s^t. 
\]
We get 
\[
v(t)=\N\left[1- \expp{-A_\sigma} \right]=\N\left[\int_0^\sigma
  \expp{-(A_\sigma -A_s)}
  dA_s\right]. 
\]
Then  we  can  replace   $  \expp{-(A_\sigma  -A_s)}$  by  its  optional
projection $B=\E^*_{\rho_s}  \left[ \expp{-A_\sigma}\right]$.  Thanks to
Lemma  \ref{lem:pm>t}, we  can  replace $A_\sigma$  in  $B$ by  $\lambda
\sigma  +  \gamma L^t_\sigma$.   Using  notations  introduced for  Lemma
\ref{lem:dlg-decomp}, we  have under $\E^*_{\mu}$ that  $\sigma=\sum_{i\in I}
\sigma^i$ (where $\sigma^i$ is the length of the excursion $\rho^i$) and
$L^t_\sigma= \sum_{i\in I} L_{\sigma^i,i}^{t-h_i}$, where $L^a_{s,i}$ is
the  local time  at level  $a$ at  time $s$  of the  exploration process
$\rho^i$.  Notice  that $H_s\leq  t$    $dA_s$-a.e.  From  Lemma
\ref{lem:dlg-decomp}, we get
\[
B=\expp{ - \int \rho_s(du) \N\left[1-\expp{ - \lambda \sigma - \gamma
      L^{t-u}_\sigma}  \right]}= \expp{- \int \rho_s(du) v(t-u)}. 
\]
Eventually, we have
\begin{align}
\nonumber
v(t)
&=\N\left[\int_0^\sigma  \expp{- \int \rho_s(du) v(t-u)}
  dA_s\right]\\
\nonumber
&=\lambda \int_0^t  da\; \expp{-\alpha_0 a} \expp{- \int_0^a dx
  \int_0^1 du \int_{(0,\infty
    )} \ell\pi(d\ell) \left[1-\expp{ - v(t-x) u \ell} \right]} \\
\nonumber
&\hspace{2cm}+(\gamma+\phi(\lambda)) \expp{-\alpha_0 t} \expp{-
  \int_0^t dx   \int_0^1 du \int_{(0,\infty
    )} \ell\pi(d\ell) \left[1-\expp{ - v(t-x) u \ell} \right]} \\
\nonumber
&=\lambda \int_0^t  da\; \expp{- \int_0^a dx
  \int_0^1 du \;\psi'(v(t-x) u) } +(\gamma+\phi(\lambda)) \expp{-
  \int_0^t dx 
  \int_0^1 du \; \psi'(v(t-x) u) } \\
\nonumber
&=\lambda \int_0^t  da\; \expp{- \int_0^a dx
  \frac{\psi(v(t-x))}{v(t-x)}} +(\gamma+\phi(\lambda)) \expp{- \int_0^t dx
  \; \frac{\psi(v(x)) }{v(x)}} \\
\label{eq:eq_v}
&=\lambda \int_0^t  da\; \expp{- \int_a^t dx
  \frac{\psi(v(x))}{v(x)}} +(\gamma+\phi (\lambda)) \expp{- \int_0^t dx
  \; \frac{\psi(v(x)) }{v(x)}},
\end{align}
where we used  formulas   \reff{eq:poisson-rh}  and \reff{eq:poisson-rl}
for the 
second equality and the convention $\frac{\psi(\infty ))}{\infty }=
\psi'(\infty )=\infty $ for the fourth. 
Notice that for all $t>0$, \reff{eq:eq_v} implies 
\[
v(t)\leq \lambda t +\gamma+\phi (\lambda).
\]
We set $v(0)=\gamma+  \phi (\lambda)$. 
Since the function $t\mapsto v(t)$ is locally bounded and measurable, we deduce 
from \reff{eq:eq_v} that $v$ is continuous and even of class $\cc^1
(\R_+)$. By differentiation w.r.t. $t$, we get 
\[
v'(t)= \lambda - \psi(v(t)).
\]
Since  $\int_{\phi (\lambda)} dv/(\lambda  -  \psi(v))=\infty $  and
$v'(0)<0$,   it  is   easy  to   check  that   $v$  is   decreasing  and
$\lim_{t\rightarrow  \infty }  v(t)  =\phi (\lambda)$.  Then  notice
that $w(t)=v(t) -\N[1-\expp{-\lambda \sigma}]= v(t) -\phi (\lambda)$
to conclude.  The case  $\lambda=0$ is similar.  The case  $\gamma=0$ is
immediate.
\end{proof}

\section{Dislocation measures}
\label{sec:dm}

For $s\in (0,\sigma)$, let  $\sigma^{s,t}$ be the size (i.e. the
Lebesgue measure) of the fragment at time $t$ which contains $s$:
\[
\sigma^{s,t}= \int_0^\sigma du\; \ind_{\{H_{[u,s]}>t\}},
\]
where $H_{[u,s]}=\min \{ H_r; r\in [s\wedge u, s\vee u]\}$. 
We consider a tagged fragment.  More precisely, let $s^*\in [0,\sigma]$ be,
conditionally  on   $\sigma$,  chosen  uniformly   and  independently  of
$\rho$.  The  process  $(\sigma^{s^*,t},  t\geq 0)$  is  a  non-increasing
process. Let $\ct_{s^*}$ the set of its jumping times. Notice there is a jump
at time $t$ if  and only if there is a node of the  CRT at level $t$ for
$s^*$,  which   is  equivalent   to  say  that   $\kappa_{s^*}(\{t\})>0$.  Let
$x(t)=(x_1(t),  \ldots)\in  \cs^\downarrow$   be  the  sequence  of  the
Lebesgue  measures   for  the  different  fragments   coming  from  the
fragmentation at time $t$ of the tagged fragment. In particular, we have
$\sum_{i=1}^\infty x_i(t) =\sigma^{s^*, t-}$.

Let $S$ be a subordinator with Laplace index $\phi $. Denote by
$(\Delta S_t, t\geq 0)$ its jumps. Let 
$\mu$  the  measure  on   $\R_+\times  \cs  ^\downarrow$  such that   for  any
non-negative measurable function, $F$, on $\R_+\times \cs ^\downarrow$,
\begin{equation}
   \label{eq:def-mu}
\int_{(0,+\infty )\times   \cs  ^\downarrow}   F(r,x)   \mu(dr,dx)=\int  \pi(dv)
  \E[F(S_v, (\Delta S_t, t\leq v))], 
\end{equation}
where $(\Delta S_t, t\leq v)$ has to be understood as the family of
jumps of the subordinator up to time $v$ ranked in non-increasing order.
Intuitively, $\mu$ is the  law of $S_T$ and the jumps of  $S$ up to time
$T$, where $T$ and $S$ are independent, and $T$ is distributed according
to  the  infinite measure  $  \pi$. Recall  $\pi_*$  is  the
``distribution''   of
$\sigma$  under $\N$  (this  is  the Lévy  measure  associated to  the
Laplace exponent $\phi $). From Theorem 9.1 in \cite{ad:falp}, we
have that $\mu(dr,dx)$ is absolutely continuous with respect to
$\pi_*(dr)$, and more precisely
 \[
r \mu(dr,dx)=\nu_r(dx)\pi_* (dr),
\]
where $(\nu_r, r>0)$ is the measurable family of dislocation measures of
the fragmentation at nodes introduced in \cite{ad:falp}.  Let $(\nu^*_r,
r>0) $ defined by $r\nu^*_r(dx)=\nu_r(dx)$ for $r>0$, so that
\begin{equation}
   \label{eq:def-nu-*}
\mu(dr,dx)=\nu^*_r(dx)\pi_* (dr).
\end{equation} 

We
refer to \cite{j:cspm} for the definition of intensity of a random point
measure.
Recall $\sigma^{s,t}  $ is the  size of the  fragment at time  $t$ which
contains $s$, and $s^*$ is uniform on $[0,\sigma]$. 
\begin{theo}
\label{theo:intensite}
The intensity  of the random  point measure
  $\sum_{t\in       \ct_{s^*}}       \delta_{t,x(t)}       (dt,dx)$       is
  given by 
  $\ind_{\{\sigma^{s^*,t-}>0\}} dt \; \nu_{\sigma^{s^*,t-}}(dx)$.
\end{theo}

\begin{rem}
\label{rem:frag-int}
  For   the    $\alpha$-stable   case,   $\psi(\lambda)=\lambda^\alpha$,
  $\alpha\in  (1,2)$, we  deduce from  Corollary 9.3  in \cite{ad:falp},
  that the fragmentation at  height is a self-similar fragmentation with
  index  $\alpha^{-1}  -1$,   and  with the  same  dislocation  measure
  $\nu_1^*$ as for  the fragmentation at node ($\nu_1$  in Corollary 9.3
  \cite{ad:falp}).  This  result   was    proved  by  Miermont
  \cite{m:sfdfstsh}, see also  \cite{m:sfdfstsn}. This result was
  previously  observed by
  Bertoin \cite{b:ssf} for the case $\alpha=2$ (in this case the exploration
  process is a reflected Brownian motion).  
\end{rem}

The rest of this section is devoted to the  proof of Theorem
\ref{theo:intensite}, which is based on  the next three lemmas. 
For a function $G$ defined on $\mathcal{S}^{\downarrow}$, $G((\Delta S_t,
t\leq r))$ has to been understood as the function $G$ evaluated on the
sequence $(\Delta S_t,
t\leq r)$ ranked  in non-increasing order and  eventually completed by
zeroes if this sequence is finite. 
\begin{lem}
\label{lem:=NgG}
Let $g\in \cb_+(\R_+)$ and $G$ a measurable non-negative function
defined on  $\mathcal{S}^{\downarrow}$.    We have 
\[
 \N\left[
\sigma \expp{-\lambda \sigma} \sum_{t\in \ct_{s^*}} g(t)G(x(t))\right]\\
=\int_{\R_+} dt\; g(t) \expp{-t 
\psi'(\phi  (\lambda))} 
\int_{(0,\infty ) \times \cs^\downarrow} \mu(dr,dx)\; r\expp{-\lambda r} G(x).
\]
\end{lem}

\begin{proof}
For $t\in\ct_{s^*}$, the quantity  $\sigma^{s^*,t-}$ is the Lebesgue measure of
$\{u\in [0,\sigma]; H_{[u,s^*]}\geq t\}$. From the property of the height
process, the set $\{u\in [0,\sigma]; H_{[u,s^*]}> t\}$ is open and can be
written as the union of $(\alpha_i, \beta_i)$, $i\in   I_t$. We
write $\sigma^{i}=\beta_i-\alpha_i$. Notice that $x(t)$ is the
sequence $(\sigma^{i}, i\in  I_t)$ ranked in non-increasing
order. Let $I_t^+$ (resp. $I_t^-$) the subset of $I_t$  of indexes such that
$\alpha_i>s^*$ (resp. $\beta_i<s^*$). Notice the sequence  $(\sigma^{i}, i\in
I_t)$ is the union of $\sigma^{s^*,t}$ and $( \sigma^{i}, i\in  I_t^+\cup
I_t^- )$.

Notice that a jump of $\sigma^{s^*,t}$ happens only is there is a node at
height $t$ in the ancestor line of $s^*$. This is equivalent to say that
$\kappa_{s^*}(\{t\})>0$. In particular, we have 
\begin{equation}
   \label{eq:sg-ig}
\sum_{t\in \ct_{s^*}} g(t)G(x(t))= \sum_{t;\kappa_{s^*}(\{t\})>0} g(t) G(
(\sigma^{i}, i\in  I_t))= \int_0^\infty
  \frac{\kappa_{s^*}(dt)}{\kappa_{s^*}(\{t\})}\;  g(t) G( 
  (\sigma^{i}, i\in  I_t)),
\end{equation}
with the convention that $\kappa(dr)/\kappa(\{r\})=0$ if $\kappa(\{r\})=0$.
We first consider 
\begin{equation}
   \label{eq:defJ}
J=\N\left[ \expp{-\lambda \sigma} \int_0^\sigma ds  \int_0^\infty
 \frac{\kappa_s(dt)}{\kappa_s(\{t\})} \; g(t)\expp{- p\sigma^{s,t}} K\Big(
 \sum_{i\in I_t^+ \cup I_t^-} \delta_{\sigma^{i}}\Big)\right],
\end{equation}
where $K$ is a non-negative measurable
function defined on the set of $\sigma$-measures on $(0,\infty )$. 
The next  computations are similar to those in Section 9.3 of
\cite{ad:falp}. 

We set $\sigma^{0,t}_-=\int_0^s du\; \ind_{\{H_{[u,s]}<t\}}$, 
$\sigma^{s,t}_-=\int_0^s du\; \ind_{\{H_{[u,s]}>t\}}$, 
$\sigma^{0,t}_+=\int_s^\sigma du\; \ind_{\{H_{[u,s]}<t\}}$, 
$\sigma^{s,t}_+=\int_s^\sigma du\; \ind_{\{H_{[u,s]}>t\}}$. Notice that 
\[
\sigma=\sigma^{0,t}_-+\sum_{i\in I_t^-} \sigma^{i} + \sigma^{s,t}_- 
+ \sigma^{s,t}_+ +\sum_{i\in I_t^+} \sigma^{i} + \sigma^{0,t}_+. 
\]
We write $K_\lambda(\mu)=K(\mu)\expp{-\lambda \langle \mu,I_d \rangle}$,
where $I_d(x)=x$ for $x\in \R_+$. 
In the integral in
\reff{eq:defJ}, we can replace 
\[
\expp{- \lambda \sigma - p\sigma^{s,t}} K\Big(
 \sum_{i\in I_t^+ \cup I_t^-} \delta_{\sigma^{i}}\Big)
=\expp{-   \lambda \sigma^{0,t}_- -
(\lambda+ p) \sigma^{s,t}_--    \lambda \sigma^{0,t}_+ -
(\lambda+ p) \sigma^{s,t}_+  } K_\lambda \Big(
 \sum_{i\in I_t^- } \delta_{\sigma^{i}}+ \sum_{i\in I_t^+ } \delta_{\sigma^{i}}\Big)
\]
 by its optional projection
\[
B= \expp{- \lambda \sigma^{0,t}_- -
(\lambda+ p) \sigma^{s,t}_-} \E^*_{\rho_s} \left[\expp{-\lambda
  \int_0^\sigma \ind_{\{H_{[0,u]}<t\}}\; du 
-(\lambda+p)
   \int_0^\sigma \ind_{\{H_{[0,u]}>t\}}\; du}  K_\lambda( 
  \sum_{i\in I_{t,+}} \delta_{\sigma^{i}}+\mu' )\right], 
\]
with   $\mu'=\sum_{i\in  I_{t,-}   }  \delta_{\sigma^{i}}   $.   Using
notations  introduced  for Lemma  \ref{lem:dlg-decomp},  we  have
$
\int_0^\sigma    \ind_{\{H_{[0,u]}<t\}}du=   \sum_{k\in    I}   \sigma_k
\ind_{\{h_k<t\}} $, $ \int_0^\sigma \ind_{\{H_{[0,u]}>t\}}du= \sum_{k\in
  I}   \sigma_k   \ind_{\{h_k>t\}}   $   and   $   \sum_{i\in   I_{t,+}}
\delta_{\sigma^{i}}= \sum_{k\in I;  h_k=t} \delta_{\sigma_k}$. Then we
deduce from Lemma \ref{lem:dlg-decomp} and \reff{eq:exp-sigma}, that
\[
B= \expp{- \lambda \sigma^{0,t}_- -
(\lambda+ p) \sigma^{s,t}_-
-\rho_s((0,t)) \phi (\lambda) -   \rho_s((t,\infty )) \phi (
\lambda+ p) }   \E[K_\lambda(\cp+\mu' )],
\]
where $\cp$ is under $\P$ a Poisson point measure  with intensity
$\rho_s( \{t\})\N[d\sigma]=\rho_s( \{t\})\pi_*(d\sigma)$. 
By    time
reversibility (see Corollary 3.1.6 in \cite{dlg:rtlpsbp}), we get
\begin{align*}
J
&=\N\Big[ \int_0^\sigma ds  \int_0^\infty
 \frac{\kappa_s(dt)}{ \kappa_s(\{t\})}\; g(t)
\expp{- \lambda \sigma^{0,t}_- -
(\lambda+ p) \sigma^{s,t}_-
-\rho_s((0,t)) \phi (\lambda) -   \rho_s((t,\infty )) \phi (
\lambda+ p) }   \\
&\hspace{5cm}  \E[K_\lambda(\cp+\mu' )]_{\big |
    \mu'=\sum_{i\in I_{t,-}} 
      \delta_{\sigma^{i}}}
\Big]\\
&=\N\Big[ \int_0^\sigma ds  \int_0^\infty
\frac{ \kappa_s(dt)}{ \kappa_s(\{t\})}\; g(t)
\expp{- \lambda \sigma^{0,t}_+ -
(\lambda+ p) \sigma^{s,t}_+
-\eta_s((0,t)) \phi (\lambda) -   \eta_s((t,\infty )) \phi (
\lambda+ p) }   \\
&\hspace{5cm}  \E[K_\lambda(\cp'+\mu' )]_{\big |
    \mu'=\sum_{i\in I_{t,+}} 
      \delta_{\sigma^{i}}}
\Big], 
\end{align*}
where $\cp'$ is under $\P$ a Poisson point measure  with intensity
$\eta_s( \{t\})\pi_*(dr)$. 
Using  similar computation as above, we get 
\[
J= \N\left[ \int_0^\sigma ds  \int_0^\infty
\frac{ \kappa_s(dt)}{\kappa_s(\{t\})} \; g(t)
\expp{-   \kappa_s((0,t)) \phi ( \lambda)-  \kappa_s((t,\infty ))
  \phi ( \lambda+ p) }   \E[K_\lambda(\cp'' )]
\right],
\]
where $\cp''$ is under $ \P$ a Poisson point measure  with intensity
$\kappa_s(\{t\}) \pi_*(d\sigma)$.
We write $f(r)$ for $ \E[K_\lambda(\cp''' )]$, where $\cp'''$ is under $
\P$ a Poisson point measure  with intensity 
$r\pi_*(d\sigma)$. Thanks to
the Poisson 
representation of \reff{eq:poisson-rh}, using 
 notation $\cn(dx,d\ell,du)=\sum_{i} \delta_{x_i, \ell_i, u_i}$, we
get 
\begin{align*}
   J
&=\E\left[ \int_0^\infty  da\expp{- \alpha_0 a  }
\sum_{x_i\leq a} \ell_i  \inv{\ell_i} g(x_i) f (\ell_i)
\expp{ - \sum_{ x_k<x_i} \ell_k \phi (\lambda)   - \sum_{a\geq
    x_j>x_i} \ell_j \phi (\lambda+ p)}   \right] \\
&=\int_0^\infty  da\expp{- \alpha_0 a  } \E\left[
\sum_{x_i\leq  a}  g(x_i) f (\ell_i)
\expp{  - x_i\int \ell\pi(d\ell) \;  
     [1-\expp{-\ell\phi (\lambda)}] -  (a-x_i)  \int \ell\pi(d\ell) \;  
     [1-\expp{-\ell\phi (\lambda+p)}]   }   \right] \\
&=\int_0^\infty  da \; \E\left[ 
\sum_{x_i\leq  a}  g(x_i) f (\ell_i) 
\expp{ - x_i \psi'(\phi (\lambda))- (a-x_i)
  \psi'(\phi (\lambda+ p)) }
 \right] \\
&= \int_0^\infty  da
\int_{(0,\infty )} \ell  \pi(d\ell) \int dx\;  \ind_{[0,a]}(x) \;  g(x)
f (\ell)  
\expp{  - x
  \psi'(\phi (\lambda)) - (a-x)
  \psi'(\phi (\lambda+p))} \\
&=\inv{\psi'(\phi (\lambda+ p))} \int_0^\infty  dt \;g(t)\expp{  - t
  \psi'(\phi (\lambda))}   \int \ell \pi(d\ell) \;  f(\ell).
\end{align*}

On the other side, let $(\Delta S_t, t\geq 0)$ be the jumps of a
subordinator $S=(S_t, t\geq 0)$ with Laplace exponent
$\phi $ and Lévy measure $\pi_* $. Standard computations yield
for $r>0$, 
\begin{align*}
  \E\left[ \expp{-\lambda S_r} \sum_{t\leq r}
  \Delta S_t 
   \expp{ - p \Delta S_t} K\big( \sum_{u\leq
   r, u\neq t} \delta_{\Delta S_u}\big)\right]
&=  \E\left[\sum_{t\leq r} \Delta S_t 
   \expp{ - (\lambda+p) \Delta S_t } K_\lambda \big(
     \sum_{u\leq    r, u\neq t } \delta_{\Delta S_u}\big) \right]\\
  & = r  \left[\int \pi_*  (d\ell) \; \ell \expp{- (\lambda+p)\ell}\right]
\E\left[{ K_\lambda\big(     \sum_{u\leq    r} \delta_{\Delta S_u}\big)  }
 \right] \\
  & = r  {\phi }' (\lambda+p)   f(r)  ,
\end{align*}
as  $\sum_{u\leq  r} \delta_{\Delta  S_u}$  is  a  Poisson measure  with
intensity    $r\pi_*(d   v)$.    Notice   that    ${\phi   }'={\phi^{-1}
}'=1/\psi'\circ \phi $ to conclude from \reff{eq:defJ} that
\begin{multline}
   \label{eq:=J=}
\N\left[ \expp{-\lambda \sigma} \int_0^\sigma ds  \int_0^\infty
 \frac{\kappa_s(dt)}{\kappa_s(\{t\})} \; g(t)\expp{- p\sigma^{s,t}} K\Big(
 \sum_{i\in I_t^+ \cup I_t^-} \delta_{\sigma^{i}}\Big)\right] \\
= \int_0^\infty  dt \;g(t)\expp{  - t
  \psi'(\phi (\lambda))}   
\int  \pi(dr)
\E\left[ \expp{-\lambda S_r} \sum_{t\leq r}
  \Delta S_t 
   \expp{ - p \Delta S_t} K\big( \sum_{u\leq
   r, u\neq t} \delta_{\Delta S_u}\big)\right].
\end{multline}
Eventually from monotone class Theorem, we get 
\begin{multline*}
    \N\left[\expp{-\lambda \sigma} \int_0^\sigma ds\; \int_0^\infty
  \frac{\kappa_s(dt)}{\kappa_s(\{t\})}\;  g(t) G( 
  (\sigma^{i}, i\in  I_t))
\right]\\
= \int_0^\infty  dt \;g(t)\expp{  - t
  \psi'(\phi (\lambda))}   
\int  \pi(dr)
\E\left[ S_r \; \expp{-\lambda S_r}  G( \Delta S_u, u\leq r\big)\right].
\end{multline*}
Then we deduce from \reff{eq:sg-ig} that
\[
   \N\left[\sigma \expp{-\lambda \sigma} \sum_{t\in \ct_{s^*}} g(t)G(x(t))\right]
= \int_0^\infty  dt \;g(t)\expp{  - t  \psi'(\phi (\lambda))}   
\int  \pi(dr)
\E\left[ S_r \; \expp{-\lambda S_r}  G( \Delta S_u, u\leq r\big)\right].
\]
where we used   that the tag  $s^*$ is chosen uniformly
on $[0,\sigma]$. We conclude by using the
definition \reff{eq:def-mu} of the measure $\mu$.
\end{proof}

\begin{lem}
\label{lem:NgH}
Let $H\in \cb_+(\R_+)$, $t>0$, $\lambda>0$.   We have 
\begin{equation}
   \label{eq:NgH}
\N\left[\sigma \expp{-\lambda \sigma}
  H(\sigma^{s^*,t-})\ind_{\{\sigma^{s^*,t-}>0\}} \right] 
=  \expp{-t
\psi'(\phi  (\lambda))} \int_{(0,\infty )}\pi_*(dr) \; r 
\expp{-\lambda r} H(r) .
\end{equation}
\end{lem}

\begin{proof}
Notice that $\sigma^{s,t-}>0$ if and only if $H_s>t$ $\N$-a.e.  Let
$p>0$.   We have 
\[
\N\left[\sigma \expp{-\lambda \sigma - p
    \sigma^{s^*,t-}}\ind_{\{\sigma^{s^*,t-}>0\}} \right] 
=\N\left[ \expp{-\lambda
    \sigma} \int_0^\sigma ds \;\expp{ -p
    \sigma^{s,t-}}\ind_{\{H_s>t\}}  \right].
\]
Similar computations as in the proof of Lemme \ref{lem:=NgG} yield
\begin{align*}
\N\left[\sigma \expp{-\lambda \sigma - p \sigma^{s^*,t-}}\ind_{\{H_s>t\}} \right]
&= \N\left[  \int_0^\sigma ds \;\expp{ -\kappa_s((0,t))
    \phi (\lambda) - \kappa_s([t,\infty ))
    \phi (\lambda+p)}\ind_{\{H_s>t\}} 
\right]  \\
&= \E\left[ \int_0^\infty  da\;\expp{- \alpha_0 a  }
\expp{ - \sum_{ x_k< t} \ell_k \phi (\lambda)   - \sum_{a\geq
    x_j\geq  t } \ell_j \phi (\lambda+ p)}   \ind_{\{a>t\}}\right] \\
&= \int_t^\infty  da
\expp{  - t  \psi'(\phi (\lambda)) -(a-t)
  \psi'(\phi (\lambda+p))} \\
&= \frac{\expp{  - t  \psi'(\phi (\lambda))}}{
  \psi'(\phi (\lambda+p))} , 
\end{align*}
where  we used  the Poisson  representation of  \reff{eq:poisson-rh} and
notation $\cn(dx,d\ell,du)=\sum_{i}  \delta_{x_i, \ell_i, u_i}$  for the
second   equality.   Since    $   \int_{(0,\infty   )}\pi_*(dr)   \;   r
\expp{-(\lambda+p)      r}      =      {\phi }'(\lambda+p)$      and
${\phi }'=1/\psi'\circ   \phi $,  the   Lemma   is  proved   for
$H(x)=\expp{-px}$ and all  $p\geq 0$. We use the  monotone class Theorem
to end the proof.
\end{proof}

\begin{lem}
\label{lem:intensite}
   We have $\pi_*(dr)$-a.e. 
\begin{equation}
   \label{eq:NgG=}
\N_r\left[
 \sum_{t\in \ct_{s^*}}
g(t)G(x(t))\right]=\N_r\left[  \int_{\R_+}
dt\; g(t)  \ind_{\{ \sigma^{s^*,t-}>0\}} \int \nu^*_{\sigma^{s^*,t-}} (dx)
\; G(x) \right].  
\end{equation}
\end{lem}

\begin{proof}
  As a  direct consequence of Lemma  \ref{lem:=NgG}, Lemma \ref{lem:NgH}
  with $H(r)=  \int_{\cs^\downarrow} \nu^*_r(dx) G(x)$ and
  \reff{eq:def-nu-*}, we have for $\lambda>0$,
\[
\N\left[\sigma \expp{-\lambda\sigma} 
 \sum_{t\in \ct_{s^*}}
g(t)G(x(t))\right]   
=\N\left[ \sigma \expp{-\lambda\sigma}  \int_{\R_+}
dt\; g(t)   \ind_{\{ \sigma^{s^*,t-}>0\}}\int \nu^*_{\sigma^{s^*,t-}} (dx)
\; G(x) \right]. 
\]
As  Laplace transforms characterize measures, and since 
the distribution of $\sigma$ under $\N$ is given by $\pi_*$, we easily
get the Lemma. 
\end{proof}

From the definition of  intensity measure (see \cite{j:cspm}), Lemma
\ref{lem:intensite} readily implies Theorem \ref{theo:intensite}.

%\bibliographystyle{abbrv}
%\bibliography{/home/delmas/cermics/Bibliographie/delmas}

\begin{thebibliography}{10}

\bibitem{ad:falp}
R.~ABRAHAM and J.-F. DELMAS.
\newblock Fragmentation associated to {L}évy processes.
\newblock {\em Preprint CERMICS [2005-291]}, 2005.

\bibitem{ap:sac}
D.~ALDOUS and J.~PITMAN.
\newblock The standard additive coalescent.
\newblock {\em Ann. Probab.}, 26(4):1703--1726, 1998.

\bibitem{b:pl}
J.~BERTOIN.
\newblock {\em L\'evy processes}.
\newblock Cambridge University Press, Cambridge, 1996.

\bibitem{b:fpcbm}
J.~BERTOIN.
\newblock A fragmentation process connected to {B}rownian motion.
\newblock {\em Probab. Th. Rel. Fields}, 117:289--301, 2000.

\bibitem{b:ssf}
J.~BERTOIN.
\newblock Self-similar fragmentations.
\newblock {\em Ann. Inst. Henri Poincar\'e}, 38(3):319--340, 2000.

\bibitem{b:smssf}
J.~BERTOIN.
\newblock On small masses in self-similar fragmentations.
\newblock {\em Stoch. Process. and Appl.}, 109(1):13--22, 2004.

\bibitem{b:rfcp}
J.~BERTOIN.
\newblock {\em Random fragmentation and coagulation processes}.
\newblock To appear, 2006.

\bibitem{bgt:rv}
N.~BINGHAM, C.~GOLDIE, and J.~TEUGELS.
\newblock {\em Regular variation}.
\newblock Cambridge University Press, Cambridge, 1987.

\bibitem{dlg:rtlpsbp}
T.~DUQUESNE and J.-F. LE~GALL.
\newblock {\em Random trees, {L}évy processes and spatial branching processes},
  volume 281.
\newblock Astérisque, 2002.

\bibitem{h:rfdsf}
B.~HAAS.
\newblock Regularity of formation of dust in self-similar fragmentations.
\newblock {\em Ann. Inst. Henri Poincar\'e}, 40(4):411--438, 2004.

\bibitem{j:cspm}
J.~JACOD.
\newblock {\em Calcul stochastique et probl\`emes de martingales}, volume 714
  of {\em Lecture Notes in Mathematics}.
\newblock Springer, 1979.

\bibitem{lglj:bplplfss}
J.-F. LE~GALL and Y.~LE~JAN.
\newblock Branching processes in {L}évy processes: Laplace functionals of snake
  and superprocesses.
\newblock {\em Ann. Probab.}, 26:1407--1432, 1998.

\bibitem{lglj:bplpep}
J.-F. LE~GALL and Y.~LE~JAN.
\newblock Branching processes in {L}évy processes: The exploration process.
\newblock {\em Ann. Probab.}, 26:213--252, 1998.

\bibitem{m:sfdfstsh}
G.~MIERMONT.
\newblock Self-similar fragmentations derived from the stable tree {I}:
  splitting at heights.
\newblock {\em Probab. Th. Rel. Fields}, 127(3):423--454, 2003.

\bibitem{m:sfdfstsn}
G.~MIERMONT.
\newblock Self-similar fragmentations derived from the stable tree {II}:
  splitting at nodes.
\newblock {\em Probab. Th. Rel. Fields}, 131(3):341--375, 2005.

\end{thebibliography}
\newcommand{\sortnoop}[1]{}

\end{document}